
\input amstex.tex
\documentstyle{amsppt}
\def\DJ{\leavevmode\setbox0=\hbox{D}\kern0pt\rlap
{\kern.04em\raise.188\ht0\hbox{-}}D}
\footline={\hss{\vbox to 2cm{\vfil\hbox{\rm\folio}}}\hss}
\nopagenumbers 
\font\ff=cmr8
\def\txt#1{{\textstyle{#1}}}
\baselineskip=13pt
\def\hf{{\textstyle{1\over2}}}
\def\a{\alpha}\def\b{\beta}
\def\B{{\bar\beta}}
\def\d{{\,\roman d}}
\def\e{\varepsilon}
\def\f{\varphi}
\def\g{\gamma} 
\def\bg{{\bar \gamma}}
\def\G{\Gamma}

\def\s{\sigma}

\def\={\;=\;}

\def\zt{\zeta(\hf+it)}

\def\D{\Delta}  
\def\no{\noindent}  
\def\R{\Re{\roman e}\,} \def\I{\Im{\roman m}\,} \def\s{\sigma}
\def\z{\zeta} 
\def\no{\noindent} 
\def\e{\varepsilon}
\def\D{\Delta} 
\def\no{\noindent} 
\def\e{\varepsilon}
\def\l{\lambda}
\def\no{\noindent} 
\font\teneufm=eufm10
\font\seveneufm=eufm7
\font\fiveeufm=eufm5
\newfam\eufmfam
\textfont\eufmfam=\teneufm
\scriptfont\eufmfam=\seveneufm
\scriptscriptfont\eufmfam=\fiveeufm
\def\mathfrak#1{{\fam\eufmfam\relax#1}}

\font\tenmsb=msbm10
\font\sevenmsb=msbm7
\font\fivemsb=msbm5
\newfam\msbfam
\textfont\msbfam=\tenmsb
\scriptfont\msbfam=\sevenmsb
\scriptscriptfont\msbfam=\fivemsb
\def\Bbb#1{{\fam\msbfam #1}}

\def \NN {\Bbb N}

\def\rightheadline{{\hfil{\ff
Sums of zeta squares}\hfil\tenrm\folio}}

\def\leftheadline{{\tenrm\folio\hfil{\ff
Aleksandar Ivi\'c }\hfil}}
\def\emptyheadline{\hfil}
\headline{\ifnum\pageno=1 \emptyheadline\else
\ifodd\pageno \rightheadline \else \leftheadline\fi\fi}

\topmatter
\title ON SUMS OF SQUARES OF THE RIEMANN ZETA-FUNCTION 
ON THE CRITICAL LINE \endtitle
\author   Aleksandar Ivi\'c \endauthor
\address{
Aleksandar Ivi\'c, Katedra Matematike RGF-a
Universiteta u Beogradu, \DJ u\v sina 7, 11000 Beograd,
Serbia (Yugoslavia).}
\endaddress
\email {\tt aivic\@matf.bg.ac.yu, 
aivic\@rgf.bg.ac.yu} \endemail
\abstract
A discussion involving the evaluation of the sum
$\sum_{0<\g\le T}|\z(\hf+i\g)|^2$ is presented, where $\g$ denotes
imaginary parts of complex zeros of $\z(s)$. Three theorems
involving certain integrals related to this sum are proved, and the 
sum is unconditionally shown to be $\ll T\log^2T\log\log T$. 
\endabstract
\endtopmatter

\head 1. Introduction
\endhead
The aim of this note is to discuss the evaluation of the sum
$$
\sum_{0<\g\le T}|\z(\hf+i\g)|^2 \eqno(1.1)
$$
and some related integrals,
where $\g$ denotes imaginary parts of complex zeros
of $\z(s)$, and where every zero is counted with its multiplicity (see
also [5] and [7]).  The interest is in obtaining unconditional
bounds for the above sum, since assuming the Riemann Hypothesis (RH)
the sum trivially vanishes.

\smallskip
A more general sum than the one in (1.1) was treated by S.M. Gonek [3]. 
He proved, under the RH, that 
$$
\sum_{0<\g\le T}\left|\z\left({1\over2}+i\left(\g + {\a\over L}\right)\right)
\right|^2
= \left(1 - \left({\sin\pi\a\over\pi\a}\right)^2\right){T\over2\pi}
\log^2T + O(T\log T)\eqno(1.2)
$$
holds uniformly for $|\a| \le \hf L$, where $L = {1\over2\pi}
\log({T\over2\pi})$. It would be interesting to recover this result
unconditionally, but our method of proof  does not
seem capable of achieving this. To evaluate the sum in (1.1) we begin
by considering the Stieltjes integral
$$
I(T) \;:=\; \int_{T_0}^{T}|\zt|^2\d S(t),\eqno(1.3)
$$
where $T_0$ is a suitable positive constant, $T > T_0$ and clearly
it may be assumed that both
$T_0$ and $T$ are not an ordinate of a zeta-zero. 
As usual (see [2, Chapter 15] or [9, Section 9.3])
$$ 
S(t) = {1\over\pi}\arg\z(\hf+it) = {1\over\pi}\I\left\{
\log\z(\hf+it)\right\} \ll \log t. \eqno(1.4)
$$
Here for $t\not=\g$ the argument of $\z(\hf+it)$ is
obtained by continuous variation along the straight lines joining the
points 2, $2+it$, $\hf+it$, starting with the value 0. If $t$ is  an 
ordinate of a zeta-zero, then we define $S(t) = S(t+0)$. As is customary 
$$
N(T) \= \sum_{0<\g\le T}1
$$
counts (with multiplicities) the number of positive imaginary parts of 
all complex zeros
which do not exceed $T$. We have (see [4, Chapter 1] or [10, Chapter 9])
$$\eqalign{
N(T) &\,= \sum_{0<\g\le T}1 = {1\over\pi}\vartheta(T) + 1 + S(T),\cr
\vartheta(T) &:= \I\left\{\log\Gamma({\txt{1\over4}}+\hf iT)\right\}
- \hf T\log\pi,
\cr}\eqno(1.5)
$$
whence $\vartheta(T)$ is continuously differentiable. In fact, by
using Stirling's formula for the gamma-function it is found that
$$
\vartheta(T) = {T\over2}\log{T\over2\pi} - {T\over2} - {\pi\over8}
+ O\left({1\over T}\right).
$$
Thus when $t \not= \g$ we can differentiate $S(t)$ by using (1.4). If $t=\g$,
then by (1.5) it is seen that $S(t)$ has a jump discontinuity which counts the 
number of zeros $\rho$ with $\g = \I\rho = t$. Let  ${\Cal J_1}(T,\e)$
denote the union of all subintervals $(\bar\g-\e,\,\bar\g+\e)\,$ lying in
$[T_0,\,T]\,$ such that $\bar\g$ denotes distinct ordinates of 
zeta-zeros, and $\e > 0$ is so small that these intervals 
are disjoint, and let
$$
{\Cal J_2}(T,\e) \:=\;[T_0,\,T] \,\backslash\,{\Cal J_1}(T,\e).
$$
Then from (1.4) and (1.5) we infer that
$$\eqalign{
I(T) &
= \lim_{\e\to0}\int_{{\Cal J_1}(T,\e)}|\zt|^2\d\left(
N(t) -{1\over\pi}\vartheta(t) - 1\right) \cr&
+ \lim_{\e\to0}{1\over\pi}\int_{{\Cal J_2}(T,\e)}|\zt|^2\I\d\log\zt\cr&
= \sum_{T_0<\g\le T}|\z(\hf+i\g)|^2 +
\lim_{\e\to0}{1\over\pi}\int_{{\Cal J_2}(T,\e)}|\zt|^2\I
 \left(i{\zeta'(\hf+it)\over\zt}\right)\d t\cr&
= \sum_{T_0<\g\le T}|\z(\hf+i\g)|^2
+ {1\over\pi}\int_{T_0}^{T}|\zt|^2\I\left(i{\zeta'(\hf+it)
\over\zt}\right)\d t,\cr}
$$
since in the last integral the zeros of $\zt$ in the denominator are 
cancelled by the zeros of $|\zt|^2 = \zt\z(\hf-it)$, and the integral 
in question in fact equals
$$
\R\int_{T_0}^{T}\zeta'(\hf+it)\z(\hf-it)\d t.
$$
Therefore we obtain the basic formula
$$\eqalign{
&\sum_{T_0<\g\le T}|\z(\hf+i\g)|^2 \cr&
= \int_{T_0}^{T}|\zt|^2\left( \d S(t) - {1\over\pi}
\I \left(i{\zeta'(\hf+it)\over\zt}\right)\d t\right).\cr}\eqno(1.6)
$$
This formula depends implicitly on the functional equation for $\z(s)$
(e.g., see [4, Chapter 1]), namely
$$
\pi^{-{s\over2}}\G(\hf s)\z(s) \= \pi^{-({1-s\over2})}\G(\hf(1-s))\z(1-s).
$$
To see this note that we have, on using (1.5),
$$\eqalign{
\sum_{T_0<\g\le T}|\z(\hf+i\g)|^2 
&= \int_{T_0}^{T}|\zt|^2\d N(t)\cr&
= \int_{T_0}^{T}|\zt|^2\left({1\over\pi}\I\left\{\hf i\,{
\Gamma'({\txt{1\over4}}+\hf it)\over 
\Gamma({\txt{1\over4}}+\hf it)}\right\}\d t
+ \d S(t)\right)\cr&
- {1\over\pi}\int_{T_0}^{T}|\zt|^2{\log\pi\over2}\d t\cr&
= {1\over\pi}\int_{T_0}^T|\zt|^2\R\{J(t)\}\d t + 
\int_{T_0}^T|\zt|^2\d S(t)\cr&
- {1\over\pi}\R\left(\int_{T_0}^T|\zt|^2{\z'(\hf+it)\over\zt}\d t\right),
\cr}
$$
say, where
$$
J(t) :=  
\hf{\Gamma'({\txt{1\over4}}+\hf it)\over \Gamma({\txt{1\over4}}+\hf it)}
+ \,{\zeta'(\hf+it)\over\zt} - \hf\log\pi.\eqno(1.7)
$$
In the functional equation we set $s = \hf + it$. 
Logarithmic differentiation gives then
$$
- \hf\log\pi + 
\hf{\Gamma'({\txt{1\over4}}+\hf it)\over \Gamma({\txt{1\over4}}+\hf it)}
+ \,{\zeta'(\hf+it)\over\zt} =
 \hf\log\pi - 
\hf{\Gamma'({\txt{1\over4}}-\hf it)\over \Gamma({\txt{1\over4}}-\hf it)}
- \,{\zeta'(\hf-it)\over\zt}.
$$
This implies that $\overline{J(t)} = -J(t)$. Hence
$J(t)$ is purely imaginary and $\R J(t) = 0$, which yields another
proof of (1.6).

\smallskip

Clearly the problem of the estimation of the sum in (1.1) reduces, on
integration by parts of the integral in (1.6) with $S(t)$, to the
estimation of the integral
$$
\int_{T_0}^T Z(t)Z'(t)S(t)\d t.
$$
Namely we have
$$
|\zt|^2 = Z^2(t),\quad Z(t) := \chi^{-{1\over2}}(\hf+it)\zt,
\quad\chi(s) = {\z(s)\over\z(1-s)},
$$
where $Z(t) \,(\in C^\infty)\,$ is so-called Hardy's function. 
The integral above, without the oscillating factor $S(t)$ is trivial, namely
it equals
$$
\hf \int_{T_0}^T (|\zt|^2)'\d t \,\ll\, T^{1/3}
$$
on using the crude bound $\zt \ll t^{1/6}$. It turns out that this
integral, with the oscillating factor $S(t)$, is more difficult.
Unconditionally the author [5] has proved
$$
\int_{T_0}^T Z(t)Z'(t)S(t)\d t \ll_\e T\log^2T(\log\log T)^{3/2+\e},
\eqno(1.8)
$$
while K. Ramachandra [7] used a different method to obtain a result
which easily implies the right-hand side of (1.8) with $T\log^2T\log\log T$. 
The same bound holds then for the sum in (1.1). In this text
we shall obtain another proof of Ramachandra's result by a method
that can be used to estimate a variety of integrals involving $S(t)$. 
We shall prove

\bigskip
THEOREM 1. {\it If $C_0$ is Euler's constant, then unconditionally}
$$\eqalign{&
\R\int_{T_0}^T \z'(\hf + it)\z(\hf - it)\d t\cr&
= -{T\over2}\log^2\left({T\over2\pi}\right) + 
(1-C_0)\log\left({T\over2\pi}\right)  +(C_0-1)T + O(T^{1/3}).\cr}
\eqno(1.9)
$$

\bigskip
THEOREM 2. {\it If the Riemann Hypothesis is true, then}
$$\eqalign{&
\int_{T_0}^T Z(t)Z'(t)S(t)\d t\cr&
= {T\over4\pi}\log^2\left({T\over2\pi}\right) + 
{C_0-1\over2\pi}T\log\left({T\over2\pi}\right) + {1-C_0\over2\pi}T 
+ O(T^{1/3}).\cr}
\eqno(1.10)
$$

\bigskip
THEOREM 3. {\it Unconditionally we have}
$$
\sum_{0<\g\le T}|\z(\hf+i\g)|^2 \;\ll \;T\log^2T\log\log T.\eqno(1.11)
$$

\bigskip
THEOREM 4. {\it Unconditionally we have}
$$\eqalign{&
\int_0^T|\zt|^2S(t)\d t \ll T\log T\log\log T,\cr&
\int_0^T|\zt|^2S^2(t)\d t \ll T\log T(\log\log T)^2,\cr}
\eqno(1.12)
$$
{\it while under the Riemann Hypothesis}
$$
\int_0^T|\zt|^2S(t)\d t \ll T\log T.\eqno(1.13)
$$

\head 2. Use of the explicit formula for $S(T)$
\endhead
A natural method to try to evaluate the  integral in (1.6) is to use 
Lemmas 1-3 of Bombieri-Hejhal [1]. This has the advantage over
other similar explicit expressions involving $S(T)$ since it incorporates
smooth functions, which is particularly satisfactory in handling the error
terms. It will be used in proving (1.13), but in evaluating the integral
(1.6) (i.e., the sum in (1.11)) this approach seems ineffective. 
Specified to $\z(s)$ ($\D = 2$, eq. (5.4) with the O(1)-term written
explicitly) [1] gives, for $|t| \ge 2$,
$t \not = \g$, 
$$\eqalign{
-{\zeta'(\s+it)\over\z(\s+it)} &=
\sum_{n=2}^\infty\Lambda(n)n^{-\s-it}v\left({\roman e}^{\log n/\log X}\right)
\cr&+ \sum_{\bar\rho}{{\tilde u}(1+({\bar\rho}-\s-it)\log X)\over 
{\bar\rho}-\s-it} + {{\tilde u}(1+(1-\s-it)\log X)\over \s+it-1},
\cr}\eqno(2.1)
$$
$$\eqalign{
&\log\zt =
\sum_{n=2}^\infty{\Lambda(n)\over\log n}
n^{-{1\over2}-it}v\left({\roman e}^{\log n/\log X}\right)  \cr&
+ \sum_{\bar\rho}\int_{1\over2}^\infty
{{\tilde u}(1+({\bar\rho}-\s-it)\log X)\over {\bar\rho}-\s-it}\d\s
+ \int_{1\over2}^\infty
{{\tilde u}(1+(1-\s-it)\log X)\over \s+it-1}\d\s.
\cr}\eqno(2.2)
$$
In (2.1) and (2.2)  ${\bar\rho}$ denotes all zeros of $\z(s)$ (complex zeros
denoted by $\rho$ and real, or ``trivial zeros" at $-2,-4,\ldots\,$),
while $X$ is a parameter satisfying $2 \le X \le t^2$
and $\Lambda(n)$ is the von Mangoldt function. We have
$$
\tilde u(s) = \int_0^\infty u(x)x^{s-1}\d x,\quad v(x) =
\int_x^\infty u(t)\d t,
$$
and $ v(0) =1$ with proper normalization,
where $u(x) \in C^\infty$ is a real positive function with compact
support in $[1,\,{\roman e}]$. One has, for every integer $k\ge0$,
$$
|\tilde u(s)| \le \max_x |u^{(k)}(x)|{\roman e}^{\max(\R s,0)+4k}
(1+|s|)^{-k},\eqno(2.3)
$$
and for every fixed integer $K>0$, and complex zeros $\rho
= \b+i\g$ and $t \ge 2$, 
$$
\eqalign{&
\sum_\rho\int_{1\over2}^\infty
{{\tilde u}(1+(\rho-\s-it)\log X)\over \rho-\s-it}\d\s\cr& \ll
1 + \sum_{|\g-t|\le1/\log X}\log\left(1 + {1\over|\g-t|\log X}\right)
+ \sum_\rho{X^{\max(\b-{1\over2},0)}\over(1+|\g-t|\log X)^K},\cr}\eqno(2.4)
$$
and for $2 \le X \le T^{3/8},\,T\ge 2,\, K \ge 3$ (for $\z(s)$ 
one can take $a = 1 -\e$ in [1, Lemma 3] by a zero-density
result of M. Jutila [6] near the line $\s=\hf$)
$$
\int_T^{2T}\sum_\rho{X^{\max(\b-{1\over2},0)}\over(1+|\g-t|\log X)^K}
\d t \ll T{\log T\over\log X},\eqno(2.5)
$$
$$
\int_T^{2T}
\left(1 + \sum_{|\g-t|\le1/\log X}\log\left(1 + {1\over|\g-t|\log X}\right)
\right)\d t\ll T{\log T\over\log X}.\eqno(2.6)
$$

\bigskip\no
Although the sum on the left-hand side of (2.4) is undefined when $t = \g$,
its absolute value is majorized, by (2.5) and (2.6),
by an integrable expression. This 
enables us to deal effectively with integrals containing $S(t)$.

\smallskip
At this point we specify $X = 3,\, T_0 = 20\, $, 
noting  that (2.4)-(2.6) will hold for $t \ge T_0$ and $T \ge T_0$. 
In (1.6) we integrate the portion with $\d S(t)$ by parts, using (1.4)
and (2.2). The integrated terms are then well defined, since $T_0 \not= \g$,
(i.e., $\z(\s+20i) \not=0$), $\,T \not= \g$. Therefore
$$
\int_{T_0}^T|\zt|^2\d S(t) = {1\over\pi}|\zt|^2\I f(t)\Big|_{T_0}^T
- {2\over\pi}\int_{T_0}^T Z(t)Z'(t)\I f(t)\d t,
$$
where $f(t)$ denotes the right-hand side of (2.2), since the integral on the 
right hand-side of the above expression exists because of (2.4)-(2.6). 
In fact, the integral in question  is $\ll T^{4/3}\log T$,
since (see [4]) $Z(t) \ll t^{1/6}$, $Z'(t) \ll t^{1/6}$.
We now integrate back by parts the above expression and obtain from (1.6)
$$
\sum_{T_0<\g\le T}|\z(\hf+i\g)|^2 
= {1\over\pi}\I\left\{\int_{T_0}^{T}|\zt|^2\left(\d f(t) - 
i{\zeta'(\hf+it)\over\zt}\d t\right)\right\}.\eqno(2.7)
$$
In (2.7)  we substitute (2.1) (with $\s = \hf$ for $\z'/\z$) and (2.2). 
Note that the terms coming from the series with $\Lambda(n)$,
the term with $\s+it-1$ and the sums over trivial zeros $-2,-4,\ldots\,$,
being continuously differentiable and well defined for any 
$t\in [T_0,\,T]\,$, will cancel out. This follows on using
$$
{\partial\f(\s+it)\over\partial t} \= i
{\partial\f(\s+it)\over\partial \s},\eqno(2.8)
$$
which holds for any holomorphic function $\f$.
The remaining terms in (2.7), which have discontinuities at $t = \g$, come 
from the sums over complex zeros $\rho$. Hence from (2.7) it follows that
$$
\sum_{T_0<\g\le T}|\z(\hf+i\g)|^2 = {1\over\pi}
\left(I_1(T) + I_2(T)\right),\eqno(2.9)
$$
where  
$$\eqalign{
I_1(T) &:= \int_{T_0}^{T}|\zt|^2\d g(t),\cr
g(t) &:= \I\sum_\rho\int_{1\over2}^\infty
{{\tilde u}(1+(\rho-\s-it)\log X)\over \rho-\s-it}\d\s\cr}
\eqno(2.10)
$$
for $t\not=\g$ and $t \ge T_0$, and in general (see (1.4) and (2.2))
for $t \ge T_0$
$$\eqalign{
&g(t) = \pi S(t) - \I\Biggl\{
\sum_{n=2}^\infty{\Lambda(n)\over\log n}
n^{-{1\over2}-it}v\left({\roman e}^{\log n/\log X}\right)  \cr&
- \sum_{k=1}^\infty  \int_{1\over2}^\infty
{{\tilde u}(1-(2k+\s+it)\log X)\over 2k+\s+it}\d\s
+ \int_{1\over2}^\infty
{{\tilde u}(1+(1-\s-it)\log X)\over \s-it - 1}\d\s\Biggr\}.
\cr}\eqno(2.11)
$$
We also have
$$\eqalign{
I_2(T) &:= \int_{T_0}^{T}|\zt|^2 h(t)\d t,\cr
h(t) &:= \I\left(i\sum_\rho
{{\tilde u}(1+(\rho-\hf-it)\log X)\over \rho-\hf-it}\right).\cr}\eqno(2.12)
$$
The integral  $I_2(T)$ converges absolutely, since the poles
of the sum over $\rho$, namely the zeros $\rho = \hf +it$, are
cancelled by the corresponding zeros of $\zt$. The integral  $I_1(T)$   in
(2.10) is an improper Stieltjes integral, which will be handled 
in the following way. An integration by parts yields
$$
I_1(T) = |\zt|^2g(t)\Big|_{T_0}^{T}
-2\int_{T_0}^{T}Z(t)Z'(t)g(t)\d t.\eqno(2.13)
$$
In the  integrated terms we can use the expression (2.10) for $g(t)$, 
since $T_0$ and $T$ are not ordinates of any zeta-zero. In the
integral on the right-hand side of (2.13) $g(t)$ is given by (2.11).
However, by using (2.4)-(2.6) we see that this integral 
converges absolutely  if $g(t)$
is given by (2.10). Hence by properties of Stieltjes integrals it is 
seen that the integration by parts in (2.13) is justified, and that 
this formula holds in fact when $g(t)$ is given by (2.10).  

\smallskip

To transform $I_2(T)$ in (2.12), let ${\bar \g}_1 < \ldots < {\bar\g}_n$
denote distinct ordinates of zeta-zeros lying in $[T_0,\,T]$, and
let $\e > 0$ be so small that all intervals $({\bar \g}_j-\e,\,
{\bar \g}_j+\e)$ $\;(j = 1,\ldots,n)$ are disjoint and lie in
$[T_0,\,T]$. Then we have
$$
\eqalign{
I_2(T) &= \left(\int_{{\bar\g}_1-\e}^{{\bar\g}_1+\e} +
\ldots + \int_{{\bar\g}_n-\e}^{{\bar\g}_n+\e}\right)|\zt|^2h(t)\d t\cr&
+ \left(\int_{T_0}^{{\bar\g}_1-\e} +
\int_{{\bar\g}_1+\e}^{{\bar\g}_2-\e} + \ldots +
\int_{{\bar\g}_{n-1}+\e}^{{\bar\g}_n-\e} +
\int_{{\bar\g}_n+\e}^{T}\right)|\zt|^2h(t)\d t\cr&
= I_{21}(T,\e) + I_{22}(T,\e),\cr}
$$
say. Hence
$$
I_2(T) = \lim_{\e\to0}(I_{21}(T,\e) + I_{22}(T,\e)) = 
\lim_{\e\to0}I_{22}(T,\e),\eqno(2.14)
$$
since $|\zt|^2h(t)$ is in fact continuous on $[T_0,\,T]$. Let now,
for $t\not=\g$,
$$
k(t) := -\I\left(\sum_\rho
\int_{1\over2}^\infty{{\tilde u}(1+(\rho-\s-it)\log X)\over \rho-\s-it}
\d\s\right),
\eqno(2.15) 
$$
Then using (2.8) we have $k'(t) = h(t)$, and integrating by parts we obtain
$$
\eqalign{
I_{22}(T,\e) &=
|\z(\hf + i({\bar\g}_1-\e))|^2k({\bar\g}_1-\e) -
|\z(\hf + iT_0)|^2k(T_0) \cr&
+ \ldots + |\z(\hf + iT)|^2k(T)
- |\z(\hf + i({\bar\g}_n+\e))|^2k({\bar\g}_n+\e)\cr&
- 2\left(\int_{T_0}^{{\bar\g}_1-\e} + \int_{{\bar\g}_1+\e}^{{\bar\g}_2-\e}
+ \ldots + \int_{{\bar\g}_n+\e}^T\right)Z(t)Z'(t)k(t)\d t.\cr}
$$
Therefore from (2.9), (2.13) and (2.14) it follows, on integrating
by parts and using  $k(t) = -g(t)\;(t\not = \g)$, that
$$\eqalign{
&\pi\sum_{T_0<\g\le T}|\z(\hf+i\g)|^2 \cr&
= \lim_{\e\to0}\sum_{T_0<\bg\le T}\Bigl(|\z(\hf+i(\bg-\e))|^2k(\bg-\e)
- |\z(\hf+i(\bg+\e))|^2k(\bg+\e)\Bigr),\cr}\eqno(2.16)
$$
where $\bg$ denotes distinct ordinates of  zeros of $\z(s)$.

\medskip
If the RH holds,  then $|\zt|^2k(t)$ is continuous at $t = \bg$,
since the zero of $\hf + i\bg - \s - it$ (at $\s = \hf,\,t = \bg$) in
the denominator of $k(t)$ is cancelled by the corresponding zero
$\z(\hf+i\bg)$. Thus the limit in (2.16) is equal to zero, and so is then
the left-hand side, which is trivial anyway. If the RH is not true,
then we can say that
$$
\pi\sum_{T_0<\g\le T}|\z(\hf+i\g)|^2 
= \sum_{T_0<\bg\le T}L(\bg),\eqno(2.17)
$$
where 
$$
L(\bg) := \lim_{\e\to0}\Bigl(|\z(\hf+i(\bg-\e))|^2k(\bg-\e)
- |\z(\hf+i(\bg+\e))|^2k(\bg+\e)\Bigr).\eqno(2.18)
$$
We have
$$
L(\bg) \;:=\;\sum_{j=1}^{r(\bg)}\ell_{\b_{j,\bg}}(\bg),
$$
where $\hf < \b_{1,\bg} < \b_{2,\bg} < \ldots \b_{r,\bg} < 1,\, 
r = r(\bg)$ are the distinct abscissas of zeros with imaginary part
$\bg$ for a given $\bg$ and, for $\b$ one of the $\b_{j,\bg}$'s, we set
$$
\eqalign{
\ell_\b(\bg) &= m(\b+i\bg)\lim_{\e\to0}\Bigl(|\z(\hf+i(\bg+\e))|^2
\I\int_{1\over2}^\infty{{\tilde u}(1+(\b-\s-i\e)\log X)\over\B-\s-i\e}\d\s
\cr&- |\z(\hf+i(\bg-\e))|^2
\I\int_{1\over2}^\infty{{\tilde u}(1+(\b-\s+i\e)\log X)\over\B-\s+i\e}\d\s
\Bigr),\cr}\eqno(2.19)
$$
where  $m(\rho)$ is the multiplicity
of the zero $\rho$. 
Namely looking at the definition (2.15) of $k(t)$ it follows that in the sum
over $\rho$, only the terms with $\rho = \b_{j,\bg}+i\bg$ will be 
discontinuous at $t = \bg \pm\e$ as $\e\to0$. 
If we develop $|\z(\hf+i(\bg\pm\e))|^2$ by using
three terms in the Taylor formula and apply to each term the reasoning
that follows, we shall obtain without difficulty that
$$
\eqalign{
\ell_\b(\bg) &= m(\b+i\bg)|\z(\hf+i\bg)|^2\lim_{\e\to0}\Bigl\{
\int_{1\over2}^2\I\Bigl({{\tilde u}(1+(\b-\s-i\e)\log X)\over\b-\s-i\e}\cr&
- {{\tilde u}(1+(\b-\s+i\e)\log X)\over\b-\s+i\e}\Bigr)\d\s\Bigr\},
\cr}\eqno(2.20)
$$
since the portion for $\s\ge2$ is continuous in $\s$. In view of
${\tilde u}({\bar s}) = \overline{{\tilde u}(s)}$  it follows that
$$\eqalign{&
\ell_\b(\bg) = m(\b+i\bg)|\z(\hf+i\bg)|^2\times
\cr& \lim_{\e\to0}\int\limits_{1\over2}^2
{2\e\R{\tilde u}(1\!+\!(\b\!-\!\s\!-\!i\e)\log X)+ 
2(\b-\s)\I{\tilde u}(1+(\b-\s-i\e)\log X)\over(\b-\s)^2+\e^2}\d\s.\cr}
$$
Noting that
$$
\eqalign{
\I {\tilde u}(1+(\b-\s-i\e)\log X) &= \I{\tilde u}(1+(\b-\s)\log X)
+ O(\e) = O(\e),\cr
\R {\tilde u}(1+(\b-\s-i\e)\log X) &= {\tilde u}(1+(\b-\s)\log X) + O(\e),
\cr}
$$
we obtain
$$
\ell_\b(\bg) = 2m(\b+i\bg)|\z(\hf+i\bg)|^2\lim_{\e\to0}\e\int_{1\over2}^2
{{\tilde u}(1+(\b-\s)\log X)\over(\b-\s)^2+\e^2}\d\s.
$$
By using ${\tilde u}(1) = 1$ (this is equivalent to $v(0) =1$) and
$$\eqalign{
&\lim_{\e\to0}\e\int_{1\over2}^2{\d\s\over(\s-\b)^2+\e^2} \cr&=
\lim_{\e\to0}\left(\arctan{2-\b\over\e} -
\arctan{\hf-\b\over\e}\right) 
= \Bigg\{\aligned {\pi\over2}\qquad
&(\b = \hf),\\ \pi\qquad&(\hf < \b \le 1),\\
\endaligned\cr}
$$
it follows that
$$
\ell_\b(\bg) = \Bigg\{\aligned \pi m(\hf + i\bg)|\z(\hf + i\bg)|^2
\;(= 0)\qquad
&(\b = \hf),\\ 2\pi m(\b + i\bg)|\z(\b + i\bg)|^2\qquad&(\b > \hf).\\
\endaligned\eqno(2.21)
$$
Therefore from (2.16), (2.18) and (2.21) we obtain, since $m(\b+i\g) = m(1-\b+i\g)$,
$$\eqalign{&
\pi\sum_{T_0<\g\le T}|\z(\hf+i\g)|^2 
= \sum_{T_0<\bg\le T}L(\bg)\cr&
= \sum_{T_0<\bg\le T,\b>{1\over2}}2\pi m(\b + i\bg)|\z(\b + i\bg)|^2
= \pi\sum_{T_0<\g\le T}|\z(\hf+i\g)|^2,\cr}
$$
which is unfortunately trivial. Thus other approaches are to be sought 
if one wishes to bound non-trivially the sum in (1.1) (cf. Theorem 3).

\smallskip
The above discussion can be clearly generalized to Dirichlet
series possessing a functional equation similar to the functional
equation satisfied by $\z(s)$ (see e.g., [8]). Also it may be
remarked that, if $\g$ is the ordinate of a zero, a similar analysis
may be made by using the identity
$$
\lim_{\e\to0}(S(\g+\e)-S(\g-\e)) = \lim_{\e\to0}(N(\g+\e)-N(\g-\e)),
$$
but again we shall obtain (by the above method of proof) nothing more
than an obvious identity.

\head 3. Proof of the Theorems
\endhead

We use (1.6) to write
$$
I(T) = \sum_{T_0<\g\le T}|\z(\hf+i\g)|^2 +
{1\over\pi}\R \int_{T_0}^T \z'(\hf + it)\z(\hf + it)\d t,
$$
where $I(T)$ is defined by (1.3).
On the other hand, using (1.5) it is found that 
$$\eqalign{
I(T) &= \int_{T_0}^T |\zt|^2\d \left( N(t) - {t\over2\pi}\log{t\over2\pi}
+ {t\over2\pi} + O\left({1\over t}\right)\right)\cr&
=  \sum_{T_0<\g\le T}|\z(\hf+i\g)|^2 
- {1\over2\pi}\int_{T_0}^T |\zt|^2\log\left({t\over2\pi}\right)\d t
+ O(1),\cr}
$$
From the two expressions for $I(T)$ it follows that
$$
\R\int_{T_0}^T \z'(\hf + it)\z(\hf - it)\d t
= - \hf\int_{T_0}^T |\zt|^2\log\left({t\over2\pi}\right)\d t
+ O(1).\eqno(3.1)
$$
One can also obtain (3.1) by using that $J(t)$, defined by (1.7), is purely
imaginary. We have (see [4, Chapter 15])
$$
\int_0^T|\zt|^2 \d t = T\log\left({T\over2\pi}\right) + (2C_0-1)T
+ E(T),\quad E(T) \ll T^c\eqno(3.2)
$$
with some $c < 1/3$ (the optimal value of $c$ is not of primary
concern here). Hence differentiating (3.2) and inserting the
resulting expression in (3.1) it is seen that the left-hand side of
(3.1) is equal to
$$
\eqalign{&
 - \hf\int_{T_0}^T \left(\log^2\left({t\over2\pi}\right) +
2C_0\log\left({t\over2\pi}\right)\right)\d t + O(T^{1/3})\cr&
= -\pi\int_{T_0/2\pi}^{T/2\pi}(\log^2u + 2C_0\log u)\d u + O(T^{1/3}).\cr}
$$
But as
$$
\int \log u\d u = u(\log u -1),\quad
\int \log^2 u\d u = u(\log^2u - 2\log u + 2),
$$
it follows that
$$
\eqalign{&
-\pi\int_{T_0/2\pi}^{T/2\pi}(\log^2u + 2C_0\log u)\d u\cr&
= -\pi\left\{
{T\over2\pi}\log^2\left({T\over2\pi}\right) + 
(2C_0-2){T\over2\pi}\log\left({T\over2\pi}\right) + 
(2-2C_0){T\over2\pi}\right\} + O(1)\cr&
= - \left\{{T\over2}\log^2\left({T\over2\pi}\right)
+ (C_0-1)T\log\left({T\over2\pi}\right) + (1-C_0)T\right\} + O(1),\cr}
$$
and (1.9) easily follows from (3.1). Note that a direct proof of (1.9),
by using approximate functional equations for $\zt$ and $\z'(\hf+it)$, seems 
quite difficult.

\medskip
To obtain Theorem 2 note that, on the RH, the left-hand side of (1.6)
vanishes. Integration by parts gives then, for any given $\e > 0$,
$$
\int_{T_0}^T Z(t)Z'(t)S(t)\d t = -{1\over2\pi}
\R\int_{T_0}^T \z'(\hf + it)\z(\hf - it)\d t + O_\e(T^\e),\eqno(3.3)
$$
since both $|Z(t)| = |\zt|$ and $Z'(t)$ are $\ll_\e t^\e$ on the RH.
Hence (1.10) follows from (1.9) and (3.3).

\medskip
To prove Theorem 3 note first that it suffices to prove that
$$
J(T) := \int_{T}^{2T}|\zt|^2\d S(t) \ll T\log^2T\log\log T.
\eqno(3.4)
$$
We have
$$ \eqalign{
J(T) &= -{1\over\pi}\int_{T}^{2T}|\zt|^2\sum_{p\le T^\delta}
p^{-1/2}\log p\cdot\cos(t\log p)\d t\cr&
+  \int_{T}^{2T}|\zt|^2\d R(t) \cr&
= -{1\over\pi}\int_{T}^{2T}|\zt|^2\sum_{p\le T^\delta}
p^{-1/2}\log p\cdot\cos(t\log p)\d t\cr&
+ i\int_{T}^{2T}\left(\zt\zeta'(\hf - it) -
\zeta'(\hf + it)\zeta(\hf - it)\right)R(t)\d t
\cr&
+ O(T^{1/3}),\cr}\eqno(3.5)
$$
where
$$
R(t) := S(t) + {1\over\pi}\sum_{p\le T^\delta}p^{-1/2}\sin(t\log p).
\eqno(3.6)
$$
Here $\delta = 1/(40k), k\in\NN $, and $p$ denotes primes. 
Note that by K.-M. Tsang's result [9] we have, uniformly in $k$,
$$
\int_T^{2T} R^{2k}(t)\d t \,\ll\, T(ck)^{2k}\eqno(3.7)
$$
with some absolute constant $c>0$. Let $V = V(T) \to\infty$ as $T \to\infty$
be a positive function, and
$$
H(T,V) := \left\{\,t \,:\, (T\le t \le 2T) \,\wedge\, |R(t)| \ge V\right\}.
$$
Then by taking $k = [V/(10c)]$ it follows from (3.7) that ($\mu(\cdot)$
denotes measure)
$$
\mu(H(T,V)) \,\ll\, T\exp(-c_1V)\qquad(c_1 > 0).\eqno(3.8)
$$
Take now $V = {100\over c_1}\log\log T$. Then by using (3.8) and
H\"older's inequality we obtain
$$
\eqalign{&
\int_{H(T,V)}(\zt\zeta'(\hf - it) - \zeta'(\hf + it)\zeta(\hf - it))R(t)\d t  
\cr&
\ll {\left(\mu(H(T,V))\int\limits_{T}^{2T}  |\zt|^4\d t
\int\limits_{T}^{2T}|\zeta'(\hf+it)|^4\d t
\int\limits_{T}^{2T}R^4(t)\d t\right)}^{1/4} \ll T,
\cr}
$$
since $\int_0^T|\zeta'(\hf+it)|^4\d t \ll T\log^6T$.
Next we have to estimate
$$
\int_{T}^{2T}|\zt|^2\sum_{p\le T^\delta}
p^{-1/2}\log p\cdot\cos(t\log p)\d t,
$$
of which the relevant part is, on using the approximate
functional equation for $\zeta^2(s)$ (see [4, Chapter 4]),
$$
\sum_{n\le 2T}d(n)n^{-1/2}\sum_{p\le T^\delta}p^{-1/2}\log p
\int_T^{2T}\cos\Bigl(t\log{t\over2\pi n}-t-{\pi\over4}\Bigr)
\sin(t\log p)\d t.
$$
If we write the trigonometric functions as exponentials, then (by the
first derivative test, namely [4, Lemma 2.1]) the above expression 
will be $O(T)$ plus two conjugate expressions, one of which is
$$
\sum_{n\le 2T}d(n)n^{-1/2}\sum_{p\le T^\delta}p^{-1/2}\log p
\int_T^{2T}\exp\Bigl(it\log{t\over2\pi np} -t\Bigr)\d t.\eqno(3.9)
$$
For $pn < C_1T$ or $pn > C_2T$ with $C_1$ sufficiently small and
$C_2$ sufficiently large the contribution will be $O(T)$ by
the first derivative test. For $C_1T \le pn \le C_2T$ we obtain,
by the second derivative test (see [4]) that the contribution is
$$
\ll \sum_{n\le 2T}d(n)n^{-1/2}\left({T\over n}\right)^{1/2}\cdot T^{1/2}
\ll T\log^2T.
$$
Finally there remains
$$\eqalign{&
\int_{[T,2T]\backslash H(T,V)}
(\zt\zeta'(\hf-it) - \z'(\hf + it)\z(\hf - it))R(t)\d t\cr&
\ll \log\log T\int_T^{2T}|\zt\zeta'(\hf-it)|\d t\cr&
\ll \log\log T{\left(\int_T^{2T}|\zt|^2\d t
\int_T^{2T}|\zeta'(\hf+it)|^2\d t\right)}^{1/2}\cr&
\ll T\log^2T\log\log T,\cr}
$$
so that the proof of Theorem 3 is complete. 

\medskip
To prove the first part of Theorem 4, namely the bound (1.12)
unconditionally, it suffices to use the Cauchy-Schwarz 
inequality for integrals and the preceding method
the proof. Then (1.12) reduces to proving that
$$
\int_{T}^{2T}|\zt|^2 \left|\sum_{p\le T^\delta}
p^{-1/2-it}\right|^2 \d t \ll T\log T(\log\log T)^2.
$$
If we write $\zt$ as a sum of Dirichlet polynomials of length $\ll \sqrt{T}$
and apply the mean value theorem for Dirichlet polynomials (see e.g.,
[4, Chapter 4]) the proof reduces to showing that
$$
\sum_{m\le X}{\omega^2(m)\over m} \ll \log X(\log\log X)^2\qquad
\Bigl(\omega(m) = \sum_{p|m}1\Bigr),
$$
which  follows  by partial summation from the formula (see [4, Chapter 13])
$$
\sum_{m\le X}\omega^2(m) = X(\log\log X)^2 + O(X\log\log X).
$$

\medskip
The problem of proving, under the RH, the bound in (1.13) is more
involved. We shall use (1.3), (2.2) and (2.4) with $\b = \hf$,
and prove the bound in (1.13) for the integral over $[T,\,2T]$,
which is sufficient. The contribution of the trivial zeros in (2.2)
as well as of the second integral is easily seen to be $\ll T\log T$,
and so is also the contribution of the first term on the right-hand
side of (2.4). The contribution of the second term is
$$\eqalign{&
\int_T^{2T}|\zt|^2\sum_{|\g-t|\le1/\log X}\log\left(1+ 
{1\over|\g-t|\log X}\right)\d t\cr&
= \sum_{T-1/\log X\le\g\le2T+1/\log X}
\int_{\g-1/\log X}^{\g+1/\log X}|\zt|^2\log\left(1+
{1\over|\g-t|\log X}\right)\d t\cr&
= {1\over\log X}\sum_{T-1/\log X\le\g\le2T+1/\log X}
\int_{-1}^1\left|\z(\hf+i(\g + {u\over\log X}))\right|^2
\log\left(1+{1\over|u|}\right)\d u.
\cr}
$$
If we exchange the order of integration and summation and use (1.2)
(which is known to hold under the RH), then the above contribution
is $\ll T\log T$, since $X = T^\delta$ and
$$
\int_{-1}^1\log\left(1+{1\over|u|}\right)\d u \ll 1.
$$
A similar analysis holds for the contribution of the second sum in
(2.4), if we split it into portions $\g \le T/3, T/3 < \g \le 3T$
and $\g > 3T$.

There remains in (2.2) the contribution of
$$
\eqalign{&
\I \sum_{n=2}^\infty {\Lambda(n)\over \log n}
n^{-1/2-it}v\left({\roman e}^{\log n/\log X}
\right) = - \sum_{2\le n\le X}{\Lambda(n)\over \log n}n^{-1/2}\sin(t\log n)
\cr& + \sum_{2\le n\le X}{\Lambda(n)\over \log n}n^{-1/2}\left(1
- v\left({\roman e}^{\log n/\log X}\right)\right)\sin(t\log n).\cr}
\eqno(3.10)
$$
By construction $v(0) = 1$ and $v'(x) = -u(x) = 0$ for $0 \le x \le 1$,
so that $v(x) = 1$ for $0 \le x \le 1$, and then $v(x)$ falls off
monotonically to $v({\roman e}) = 0$. We have
$$
1  - v\left({\roman e}^{\log n/\log X}\right) \ll {\log n\over\log X}
\quad (2\le n \le X).
$$
By the analysis that follows, the contribution of the right-hand side
of (3.10) containing the $v$-function will be $\ll T\log T$, so
we may concentrate on the first sum the right-hand side of (3.10).
Since $\Lambda(n) = \log p$ if $n = p^m,\, p$ a prime and zero otherwise,
it seen by using the approximate functional equation for $\z^2(s)$
(see e.g., [4, Chapter 4])
the main contribution from this sum will be contained in a multiple of
$$
\sum_{n\le T/\pi}d(n)n^{-1/2}\sum_{p\le T^\delta}p^{-1/2}
\int_T^{2T}\cos\left(t\log{t\over2\pi n} - t - {\pi\over4}\right)
\sin(t\log p)\d t.
$$
By writing the trigonometric functions as exponentials and using the
first and second derivative test it is seen that the above expression
equals $O(T\log T)$ plus 
$$
\sum_{n\le T/\pi}d(n)n^{-1/2}\sum_{p\le T^\delta}p^{-1/2}
\I\left\{\int_T^{2T}{\roman e}^{iF(t)}
\d t\right\},\eqno(3.11)
$$
where we have set
$$
F(t) := t\log{t\over2\pi np} - t - {\pi\over4},
$$
so that
$$
F'(t) = \log{t\over2\pi np},\quad F''(t) = {1\over t}.
$$
There will be a saddle point $t_0$ (solution of $F'(t_0) = 0$) for
$t_0 = 2\pi np$. We split the range of summation in (3.11) into the
subranges: I) $2\pi np \le T - H$, II) $T - H < 2\pi np \le T + H$,
III) $T + H < 2\pi np \le 2T-H$, IV) $2T - H < 2\pi np \le 2T + H$
and V) $2\pi np > 2T + H$. The choice for $H$ will be
$$
H \= T^{2/3}.
$$
The contribution of the ranges I) and V) is estimated analogously. In the
former we have
$$
F'(t) \ge \log {T\over 2\pi np} \ge \log {T\over T - H} \sim {H\over T}
= T^{-1/3}.
$$
Hence by the first derivative test the contribution is
$$
\ll T^{1/3}\sum_{p\le T^\delta}p^{-1/2}\sum_{n\le T/\pi} d(n)n^{-1/2}
\ll T^{5/6+\delta} \ll T\log T.
$$
The contribution of the ranges II) and IV) is also estimated analogously.
In the former we estimate the integral as $\ll T^{1/2}$ by the second
derivative test. The contribution is then
$$
\eqalign{&
\ll T^{1/2}\sum_{p\le T^\delta}p^{-1/2}\sum_{{T-H\over2\pi p}\le n
\le {T+H\over2\pi p}}d(n)n^{-1/2}\cr&
\ll T^{1/2}\sum_{p\le T^\delta}p^{-1/2}(T/p)^{-1/2}Hp^{-1}\log T
\ll T^{2/3}\log T\log\log T.\cr}
$$
There remains the range III) in which the saddle point method is
used in the form of [10, Lemma 4.6]  with the first $O$-term there in the
form
$$
O(\l_2^{-1}\l_3^{1/3}) = O\left(T\cdot\left({1\over T^2}\right)^{1/3}\right)
= O(T^{1/3}).
$$
In view of the choice of $H$ we have that ($t_0 = 2\pi np$)
$$\eqalign{
\I\left\{\int_T^{2T}{\roman e}^{iF(t)}\d t\right\}
&= \I \sqrt{2\pi\over F''(t_0)}{\roman e}^{iF(t_0)+{1\over4}\pi i}
+ O(T^{1/3})\cr&
= \I \sqrt{2\pi\over F''(t_0)}{\roman e}^{-2\pi inp} + O(T^{1/3})
= O(T^{1/3}),\cr}
$$
since ${\roman e}^{-2\pi inp} = 1$ is real. This contribution is then
$$
\ll T^{1/3}\sum_{p\le T^\delta}p^{-1/2}\sum_{{T\over 2\pi p}\le n
\le {T\over \pi p}}d(n)n^{-1/2} \ll T^{5/6}\log T\log\log T.
$$
Collecting the above estimates we obtain (1.13), and the proof of
Theorem 4 is complete. 

\medskip
In the preceding proof the RH was used 
via Gonek's formula (1.2) and (2.4) with $\b = \hf$. Thus
an unconditional proof of (1.13) would require an unconditional proof
of (1.2) (or an adequate upper bound estimate) plus an estimation
of the integral with $|\zt|^2$ and (2.4) without the simplifying
condition $\b = \hf$.

\vfill
\break
\topglue1cm

\Refs

\bigskip

\item{[1]} E. Bombieri and D.A. Hejhal, {\it On the distribution of zeros
of linear combinations of Euler products}, Duke Math. J. {\bf80}(1995),
821-862.

\item{[2]} H. Davenport, {\it Multiplicative Number Theory }
(2nd edition), GTM{\bf74}, Springer, New York-Heidelberg-Berlin, 1980.

\item{[3]} S.M. Gonek, {\it Mean values of the Riemann zeta-function and
its derivatives}, Invent. math. {\bf75}(1984), 123-141.

\item{[4]} A. Ivi\'c, {\it The Riemann zeta-function}, John
Wiley \& Sons, New York, 1985.

\item{[5]} A. Ivi\'c, {\it On certain sums over ordinates of zeta-zeros},
Bulletin CXXI de l'Acad\'emie Serbe des Sciences et des 
Arts - 2001, Classe des Sciences math\'ematiques et naturelles, 
Sciences math\'ematiques No. {\bf26}, pp. 39-52.

\item{[6]} M. Jutila, {\it Zeros of the zeta-function near the
critical line}, in ``Studies in Pure Mathematics. To the memory
of Paul Tur\'an", Birkh\"auser, Basel-Stuttgart, 1982, 365-394.

\item{[7]} K. Ramachandra, {\it On a problem of Ivi\'c},
Hardy-Ramanujan Journal {\bf23}(2001), 10-19.

\item{[8]} A. Selberg, {\it  Old and new conjectures and results
about a class of Dirichlet series}, Proc. Amalfi Conf. Analytic
Number Theory, eds. E. Bombieri et al., Universit\`a di Salerno,
Salerno, 1992, 367-385.

\item{[9]} K.-M. Tsang, {\it Some $\Omega$--theorems for the Riemann
zeta-function}, Acta Arithmetica {\bf 46}(1986), 369-395.

\item{[10]} E.C. Titchmarsh, {\it The theory of the Riemann
zeta-function}, 2nd edition, Oxford University Press, Oxford, 1986.

\vskip2cm

Aleksandar Ivi\'c

Katedra Matematike RGF-a

Universitet u Beogradu

\DJ u\v sina 7, 11000 Beograd

Serbia and Montenegro

e-mail: {\tt aivic\@rgf.bg.ac.yu, eivica\@ubbg.etf.bg.ac.yu}

\endRefs

\bye